\documentclass{article}
\usepackage{mathrsfs}
\usepackage{amsfonts,amsthm,amssymb,amsmath}
\usepackage{graphicx,color,dsfont}
\usepackage{longtable}
\usepackage{rotating}
\usepackage{multirow}
\usepackage{float}
\usepackage{bm}
\usepackage{bbm}
\usepackage{enumitem}
\usepackage{theoremref}
\usepackage{algorithm}
\usepackage{algpseudocode}
\usepackage{color}
\usepackage{CJKutf8}
\usepackage{booktabs}
\usepackage{verbatim}

\usepackage[top=1in,bottom=1in,left=1in,right=1in, a4paper
]{geometry}  

\usepackage{tikz}
\usetikzlibrary{arrows,shapes}

\usepackage[normalem]{ulem}

\numberwithin{equation}{section}

\def\bx{{\bf x}}
\def\bs{{\bf s}}
\def\be{{\bf e}}
\def\bK{{\bf K}}
\def\by{{\bf y}}

\usepackage{natbib}

\usepackage{color}

\usepackage{hyperref} 
\hypersetup{
    colorlinks=true,
    linkcolor=blue,
    urlcolor=red,
    linktoc=all,
    citecolor=blue
           }

\title{Review of Large-Scale Simulation Optimization}

\author{Weiwei Fan\\ {\small School of Economics and Management, Tongji University, Shanghai, China}\\[5pt]
L. Jeff Hong\\
{\small School of Management and School of Data Science, Fudan University, Shanghai, China}\\[5pt]
Guangxin Jiang\\
{\small School of Management, Harbin Institute of Technology, Harbin, Heilongjiang Province, China}\\[5pt]
Jun Luo\\
{\small Antai College of Economics and Management, Shanghai Jiao Tong University, Shanghai, China}}
\date{\today}

\begin{document}

\maketitle

\begin{abstract}
    \noindent Large-scale simulation optimization (SO) problems encompass both large-scale ranking-and-selection problems and high-dimensional discrete or continuous SO problems, presenting significant challenges to existing SO theories and algorithms. This paper begins by providing illustrative examples that highlight the differences between large-scale SO problems and those of a more moderate scale. Subsequently, it reviews several widely employed techniques for addressing large-scale SO problems, such as divide and conquer, dimension reduction, and gradient-based algorithms. Additionally, the paper examines parallelization techniques leveraging widely accessible parallel computing environments to facilitate the resolution of large-scale SO problems.
\end{abstract}

\section{Introduction}
Simulation optimization (SO) pertains to a category of optimization problems in which the objective and/or constraints lack analytical forms and can only be evaluated through running simulation experiments based on a well-specified simulation model. Specifically, we address SO problems in the form
\begin{equation}\label{eq:so}
    \max_{\mathbf{x} \in \mathbb{X}}\,\{f(\mathbf{x}):=\mathbb{E}[F(\mathbf{x})]\},
\end{equation}
where  $\mathbf{x}$ represents the vector of decision variables; $\mathbb{X}$ denotes the feasible set, often defined by deterministic or explicit functions; and $F(\mathbf{x})$ is a real-valued random variable, reflecting the stochastic response of the simulation model evaluated at  $\mathbf{x}$. The distribution of  $F(\mathbf{x})$ is an unknown function of  $\mathbf{x}$, but its realizations may be observed through running simulation experiments at $\mathbb{X}$. There are numerous instances of SO problems. For example, in inventory management, one may determine the base-stock levels of the products to minimize the expected total inventory cost, using a simulation model that captures the dynamics of an inventory system and calculates the total (stochastic) inventory cost over a time period \citep{glasserman1995sensitivity,wang2023large}. In hospital management, one may need to determine the schedule of the surgeries to minimize the expected total waiting and overtime cost \citep{fan2020a}. 

Compared to other formulations of stochastic optimization, SO presents distinct advantages and disadvantages. The primary advantage lies in the generalizability of simulation modeling, enabling the modeling of complex process dynamics and the incorporation of intricate operational details that are otherwise challenging to represent. Consequently, SO can theoretically be applied to any problem for which a simulation model can be constructed. Another clear advantage is the segregation of modeling and optimization. In most optimization techniques, models are constructed based on the techniques used to solve them. For example, linear programming solvers necessitate the modeling of objectives and constraints using linear functions. In contrast, simulation models are often constructed for multiple purposes and can be utilized across various optimization problems. This segregation allows modelers to create more realistic models and utilize them repeatedly. However, these advantages also present challenges for SO. The generalizability of simulation modeling often results in models with lower levels of abstraction and greater incorporation of dynamics and details. Consequently, evaluating simulation models is significantly more time-consuming compared to models with explicit formulas commonly used in other stochastic optimization problems. As a result, SO problems typically require considerably more computational effort to solve. Additionally, the segregation of modeling and optimization often leads to a lack of effective structure for optimization. Consequently, simulation models are frequently treated as black boxes, and structural information beyond gradients is seldom utilized in SO algorithms.

According to the structure of the feasible set $\mathbb{X}$, SO problems of Equation (\ref{eq:so}) may be categorized as follows:
\begin{itemize}
\item Ranking and Selection (R\&S): In this category, $\mathbb{X}$ is a finite set with a relatively small number of solutions (alternatives), and all the solutions can be simulated to identify the best solution with high confidence.

\item Continuous SO: In this category, the decision variable $\mathbf{x}$ is a continuous vector, and the feasible set $\mathbb{X}$ is a subset of $\mathcal{R}^d$.

\item Discrete SO: In this category, the decision variable $\mathbf{x}$ is discrete (often integer-ordered), and the feasible set $\mathbb{X}$ is an integer lattice.
\end{itemize}
Readers may refer to \cite{hong2009brief} for a more detailed description of these three categories. It is important to note that these divisions are not exclusive. For example, there are mixed integer SO problems that include both continuous and discrete decision variables \citep{wang2012retrospective}. Furthermore, the divisions are sometimes blurred depending on the formulation as well as the computational budget/time. For example, for a problem with continuous decision variables, one may discretize the feasible region into a set of grid points. Then, the problem may be solved either as an R\&S problem or a discrete SO problem, often depending on whether there is sufficient computation budget/time to simulate all grid points.

There are numerous outstanding reviews on the problems of Equation (\ref{eq:so}). \cite{fu2002optimization} and \cite{hong2009brief} are excellent review articles that provide comprehensive overviews of the development of the entire field of SO at their respective times. Additionally, \cite{fu2015handbook} is a handbook containing many articles that provide excellent reviews on different subareas/methodologies of SO. In recent years, there has been a surge in reviews focusing on different subjects of SO due to their increasing popularity. For example, \cite{frazier2018tutorial} reviewed Bayesian optimization techniques commonly used to solve various SO problems, \cite{hong2021surrogate} reviewed surrogate-based SO methods popular in handling expensive simulations, and \cite{hong2021review} reviewed the recent development of R\&S as it becomes increasingly popular. Furthermore, there have been a few excellent reviews on the applications of SO in specific fields, such as \cite{wang2023simulation} reviewing the application of SO methods to healthcare resource planning, \cite{zhou2021classification} reviewing SO in the field of maritime logistics, and \cite{jalali2015simulation} reviewing SO in inventory replenishment. Most recently, with the rise of artificial intelligence (AI), some SO methods have been applied in AI, and \cite{peng2023simulation} reviewed how these methods underpin modern AI techniques.

In this paper, we focus on large-scale SO problems, also in the form of Equation (\ref{eq:so}). There are several reasons for this focus. Firstly, due to the rapid development of computational technologies, such as the widespread availability of massive parallel computing hardware like graphic processing units (GPUs) and the easy accessibility of efficient computational software frameworks like TensorFlow \citep{tensorflow2015-whitepaper}, we are now able to simulate more complex and larger systems, explore larger solution spaces, and consider more decision variables. All of these factors contribute to larger-scale SO problems, with either a larger number of feasible solutions, higher-dimensional decision variables, or both. Secondly, previous research has repeatedly demonstrated that large-scale problems are inherently challenging and cannot be simply extrapolated from moderate-scale problems. Algorithms that are efficient for moderate-scale problems may perform poorly for large-scale problems, and vice versa. Often, large-scale problems require a completely different mindset or framework for solution and warrant careful reconsideration of related theories and algorithms. Thirdly, the use of parallel computing is essential for solving large-scale problems. However, efficiently leveraging parallel computing for large-scale SO is a non-trivial task. A thorough review of the related literature is itself very valuable. For these reasons, this paper chooses to focus on large-scale SO problems and their solution algorithms. In Section 2, we discuss the differences between large-scale SO problems and moderate-scale problems. In Sections 3, 4, and 5, we review three types of frameworks/algorithms for solving large-scale SO problems: divide and conquer, dimension reduction, and gradient-based algorithms. The issue of efficient parallelization in SO is considered in Section 6, followed by conclusions in Section 7.

Before proceeding to the next section, it is important to note that some problems related to SO do not conform to the framework presented in Equation (\ref{eq:so}). For example, in a simulation-based decision-making scenario, one might be interested in multiple performance measures derived from a simulation model. For instance, in the previously mentioned inventory example, one might seek to consider both the inventory cost and the fill rate. Such problems can be formulated as either multi-objective SO problems (refer to \cite{hunter2019introduction} for a recent review) or constrained SO problems with stochastic constraints involving simulation output (see \cite{andradottir2010fully} for an article on constrained R\&S and \cite{hong2015chance} for an article on chance-constrained programs). Another type of SO problem not covered by the framework of Equation (\ref{eq:so}) is online SO, also known as contextual SO or SO with covariates. The primary distinctions between online SO and conventional problems lie in the ability to observe and utilize new information during problem-solving, as well as the real-time or near real-time nature of the problems. Analogous to Problem (\ref{eq:so}), the online SO problem may be formulated as:
\begin{equation}\nonumber
\max_{\mathbf{x} \in \mathbb{X(\mathcal{A})}}\{f_{\mathcal{A}}(\mathbf{x}):=\mathbb{E}[F(\mathbf{x})|\mathcal{A}]\},
\end{equation}
where $\mathcal{A}$ represents the newly observed information during problem-solving. The information set $\mathcal{A}$ can take various forms. A common scenario is when $\mathcal{A}$ is in the form of $\{\mathbf{Y} = \mathbf{y}\}$, and the objective function is given by $\mathbb{E}[F(\mathbf{x})|\mathbf{Y} = \mathbf{y}]$, where $\mathbf{Y}$ represents the covariates in the simulation model with observed values $\mathbf{y}$—i.e., the simulation model parameters that are only revealed at the moment. For such problems, \cite{hongjiang2019} propose the general framework of ``offline learning online optimization" to address them. Interested individuals may also refer to \cite{shenRSCovariates2021} and \cite{du2023} for algorithms that solve online R\&S problems.

\section{Large Is Different}

For SO problems,  the term ``scale" can have varying interpretations depending on the type of problem.  In the context of R\&S problems where all solutions are simulated, the scale parameter pertains to the number of solutions. Classical R\&S algorithms typically address problems with up to a few hundred solutions, which are referred to as moderate-scale problems in this paper. In recent R\&S literature, large-scale problems are characterized by thousands to even millions of solutions. However, it is important to note that the classification of large scale is influenced by available computing power. For example, a paper by \cite{nelson2001simple} considers R\&S problems with 500 solutions as large-scale, but with advancements in computing technology, such problems are now classified as moderate-scale.  In the context of continuous and discrete SO problems, the scale parameter often denotes the dimension of the decision vector, i.e., the number of decision variables. Problems considered moderate-scale in the literature typically involve up to ten dimensions, while large-scale problems may encompass tens to even hundreds or thousands of dimensions. For example, \cite{wang2023large} consider an extremely large-scale continuous SO problem that has 500,000 decision variables. It is important to observe that the size of the feasible region typically grows exponentially with the dimension. Consequently, SO problems with high dimensions are generally more challenging than those with large feasible regions.

For moderate-scale SO problems, the scale parameter is often not explicitly taken into account when designing algorithms, or at least is not considered as crucial as some other parameters. However, it is frequently observed that overlooking this factor may lead to unexpected consequences when applying these algorithms to solve large-scale problems. In other words, algorithms that are effective and efficient for moderate-scale problems may perform poorly for large-scale problems. Conversely, algorithms that perform poorly for moderate-scale problems may surprisingly perform well for large-scale problems. These observations convince us that large-scale SO problems are fundamentally different. They cannot be simply extrapolated or straightforwardly extended from moderate-scale problems. They require different theories, different algorithms, and most importantly, different mindsets.

In the remainder of this section, we present four examples from our research experience to illustrate the distinctions between large-scale and moderate-scale SO problems, and to emphasize the need for different mindsets to tackle these challenges.

\subsection{Sample Optimality for Fixed-Budget R\&S}\label{subsec:RS}

Consider a R\&S problem in which the feasible set is denoted as $\mathbb{X}=\{\mathbf{x}_1,\mathbf{x}_2,\ldots,\mathbf{x}_k\}$. Without loss of generality, let $\mathbf{x}_k$ represent the optimal solution of Problem (\ref{eq:so}), and the objective is to correctly select $\mathbf{x}_k$ from $\mathbb{X}$.  Suppose that we have a total budget of $N$ simulation observations, and our goal is to allocate them to the $k$ feasible solutions to maximize the probability of correct selection. This problem, known as a fixed-budget R\&S problem, has been extensively studied in the literature \citep{hong2021review}.

To simplify notation, we use $\bar F_i(n_i)$ to represent the sample mean of $F({\mathbf{x}_i})$ calculated with $n_i$ independent simulation observations at $\mathbf{x}_i$ for all $i=1,2,\ldots,k$. Many R\&S algorithms aim to approximate the optimal solution to the following problem:
\begin{eqnarray}
    &\max& \mathbb{P}\left\{ \bar F_k(n_k)\ge \bar F_i(n_i),\ \forall\ i=1,\ldots,k-1 \right\} \label{eq:ocba}\\
    &{\rm s.t.}& n_1+n_2+\cdots+n_k = N \nonumber\\
    && n_i\ge 1,\quad \forall\ i=1,\ldots,k.  \nonumber
\end{eqnarray}
Here, the objective is the probability of correct selection (PCS), as the solution with the largest sample mean is typically selected as the best. It is important to note that even when $F({\mathbf{x}_i})$ follows a normal distribution with known mean and variance for all $i=1,\ldots,k$, Problem (\ref{eq:ocba}) still has no closed-form solution. Many fixed-budget R\&S algorithms first find an approximate solution to the problem and then sequentially approximate the approximate solution by learning the unknown distributional parameters (e.g., the means and variances under the normal assumption). These algorithms include the optimal computing budget allocation (OCBA) algorithms \citep{chen2000simulation}, the large-deviation algorithms \citep{glynn2004large}, and many variants of these two influential algorithms.

While algorithms that approximate the optimal solution of Problem (\ref{eq:ocba}) are highly efficient for moderate-scale fixed-budget R\&S problems, \cite{hong2022solving} have demonstrated that even the allocation based on the true optimal solution of Problem (\ref{eq:ocba}) performs poorly for large-scale problems. They have shown that if the total budget $N$ grows at a rate slower than $k\log k$, the optimal PCS achieved by Problem (\ref{eq:ocba}) tends to zero. For example, if the total budget grows linearly in $k$ (e.g., the per solution budget is fixed), the PCS of large-scale problems (i.e., when $k$ is large) approaches zero. Numerical evidence supports these findings and indicates that algorithms based on Problem (\ref{eq:ocba}), such as the OCBA algorithm, may perform poorly when $k$ reaches $10^5$ or more \citep{hong2022solving,li2023surprising}.

\cite{hong2022solving} demonstrate that, in order to ensure that the PCS remains bounded away from zero, the total budget must increase at least in the order of $k$. They refer to algorithms that achieve this order as ``sample optimal" for large-scale fixed-budget R\&S problems. They have developed an algorithm based on the knockout tournament concept proposed by \cite{zhong2022knockout} to attain this ``sample optimality". According to this definition of sample optimality, it is evident that the OCBA algorithms and the large-deviation algorithms do not meet the sample optimality criteria. The numerical findings of \cite{hong2022solving} indicate that the new algorithm indeed exhibits significantly better performance for large-scale problems, even though its performance for moderate-scale problems is not as competitive as the OCBA algorithm.

Recently, \cite{li2023surprising} discover and prove that the greedy algorithm, which always allocates the next simulation observation to the current sample best, is sample optimal for large-scale fixed-budget R\&S problems. This is a truly surprising finding because the greedy algorithm is generally perceived as a na\"ive and inefficient algorithm for R\&S problems. However, the finding is supported by numerical evidence. Such a simple and na\"ive algorithm typically has very poor performance for moderate-scale problems, yet it may perform significantly better than the OCBA algorithm when $k$ is large.

\subsection{Uniform Sampling of COMPASS} \label{subsec:compass}

The COMPASS algorithm, as presented by \cite{hong2006discrete}, is one of the most popular and efficient random search algorithms for discrete SO problems \citep{li2015mo}. At each iteration, the algorithm establishes a most promising area (MPA) surrounding the current best sample solution. This MPA is a closed and bounded polyhedron, defined by linear constraints, encompassing all feasible solutions closer to the current best sample solution than any other previously visited solutions. The algorithm then proceeds to randomly sample (and simulate) a handful of feasible solutions from the MPA using a {\it uniform sampling distribution} over the feasible solutions within the MPA \citep{hong2006discrete}.

The rationale for employing a uniform sampling distribution lies in the fact that the MPA, as defined, encompasses all solutions that appear more promising than others. It seems natural, therefore, to uniformly sample from these promising solutions in order to discover better solutions. COMPASS works well in solving discrete SO problems of up to ten dimensions. However, \cite{hong2010speeding} note that the algorithm's performance deteriorates rapidly as dimensionality increases, with the algorithm typically exhibiting poor performance for problems featuring more than 15 dimensions.

The study by \cite{hong2010speeding} reveals an intriguing finding: the decline in performance of COMPASS for high-dimensional problems is attributed to uniform sampling, which appears innocuous at first. The MPA operates on the basic premise that better solutions are more likely to be found among those that are close to the current best sample. However, as the dimensionality increases, the proportion of solutions within the MPA that are close to the sample best diminishes rapidly, rendering it highly improbable for the uniform sampling to identify these solutions. This discovery can be illustrated by considering a $d$-dimensional sphere inscribed in a $d$-dimensional cube. The cube represents the MPA, while the sphere denotes the desirable solutions. It becomes evident that the proportion of the sphere to the cube (i.e., the ratio of the sphere's volume to the cube's volume) decreases swiftly as $d$ increases. For instance, the ratio is 0.785 when $d=2$, 0.081 when $d=6$, 0.0025 when $d=10$, $2.46\times 10^{-8}$ when $d=20$, and $1.54\times 10^{-28}$ when $d=50$. This simple model elucidates why uniform sampling is ineffective in identifying better solutions for high-dimensional problems. It is worth noting that \cite{hong2010speeding} utilized a slightly more complex model to represent the MPA, but the fundamental concept remains unchanged. To address this issue, \cite{hong2010speeding} propose replacing uniform sampling with coordinate sampling, which uniformly samples a coordinate direction and then uniformly samples from the feasible solutions within the MPA along the coordinate direction from the current best sample solution. They demonstrate that this simple modification significantly enhances the performance of the resulting COMPASS algorithm for high-dimensional discrete SO problems, enabling it to handle problems with up to 50 dimensions.

\subsection{Higher-Order Smoothness Matters}\label{subsec:smoothness}

For continuous SO problems, the smoothness of the objective function with respect to the decision variables has always played a crucial role in algorithm design. However, most algorithms only require up to the second-order smoothness (e.g., second-order continuous differentiability). For example, convex SO algorithms or local SO algorithms may require the gradient or sometimes the Hessian, to guild the search process \citep{robbins1951stochastic, KW1952stochastic, spall2009feedback} or to construct local quadratic approximations \citep{chang2013stochastic}; while global SO algorithms may need continuity or differentiability to achieve certain levels of rates of convergence \citep{wang2023gaussian, zhang2022actor}.

For high-dimensional continuous SO problems, global SO algorithms often encounter the curse of dimensionality. For example, \cite{wang2023gaussian} demonstrate that the Gaussian-process based search algorithms have a worst-case rate of convergence of approximately $N^{-1/(2+d)}$, where $N$ is the total number of simulation observations and $d$ is the dimension of the decision vector. Consequently, as the dimension $d$ increases, the rate of convergence deteriorates rapidly. This is not unexpected, as global SO algorithms typically utilize a surrogate model to guide the random search and classical nonparametric surrogate models that approximate the objective function often have an optimal worst-case rate of convergence of $N^{-c_1/(c_2+d)}$, where $c_1$ and $c_2$ are small positive constants \citep{stone1982optimal, doring2018rate, gao2020towards}. These classical nonparametric surrogate models usually require up to second-order smoothness and encounter significant challenges when fitting high-dimensional surfaces.

The recent research on high-dimensional surface fitting, such as reproducing kernel Hilbert space (RKHS), has reshaped our understanding of the curse of dimensionality. It is now evident that the dimension itself is not necessarily the primary challenge in high-dimensional surface fitting. Instead, the difficulty may be influenced by a combination of dimension and smoothness. If the smoothness parameter $\nu$ of the surface is known, numerous studies have demonstrated that certain surface fitting algorithms can achieve a rate of convergence of approximately $n^{-\nu/(2\nu+d)}$. For example, see \cite{raskutti2014early}, \cite{yang2016bayesian}, \cite{hamm2021adaptive}, and \cite{ding2023random}. This finding indicates that, when the smoothness of the target surface is comparable to the dimension $d$, the curse of dimensionality does not apply. In the extreme scenario where the target surface is infinitely smooth (i.e., $\nu=\infty$), these nonparametric surface fitting algorithms can achieve a rate of convergence of approximately $n^{-1/2}$, which aligns with the performance of many parametric surface fitting algorithms such as linear regression.

When employing these surface fitting algorithms in surrogate-based SO algorithms, a higher rate of convergence may also be achieved. However, the challenge often lies in integrating these surface fitting algorithms with the optimization algorithms. A recent study by \cite{ding2021high} demonstrated that by combining a sparse grid and kernel ridge regression with a Brownian field kernel, their SO algorithm may achieve an approximate rate of convergence of $N^{-1/6}$ or $N^{-3/10}$, under respective smoothness assumptions, which does not depend on the dimension of the decision vector. Their numerical results indicate that the algorithm performs effectively for test problems of up to 100 dimensions.

\subsection{Cost of IPA}\label{subsec:ipa}

The abbreviation ``IPA" in this context does not stand for India pale ale, but rather refers to ``infinitesimal perturbation analysis," a popular approach for computing the sample-path gradient of a stochastic simulation model \citep{ho1983infinitesimal, fu2008you}. By leveraging the sample-path gradient, one can apply popular stochastic gradient descent (SGD) algorithms to solve continuous SO problems, including those with high-dimensional decision spaces. For example, in their work, \cite{wang2023large} utilize SGD algorithms to solve a simulation-based inventory optimization problem featuring up to 500,000 decision variables.

The IPA method computes a gradient estimator, i.e., a sample-path gradient, concurrently with the simulation process. In the classical literature of simulation gradient estimation, IPA is recognized as an efficient algorithm, providing an unbiased gradient estimator with minimal computational overhead when applicable \citep{fu2008you}. However, in situations where the dimension of the decision vector is high, the resulting high-dimensional gradient may lead to a significant computational overhead, diminishing the efficiency and desirability of IPA.

To illustrate the concept, let's consider a straightforward dynamic simulation that computes the stochastic function $F(\mathbf{x})$, which can be represented as follows:
\begin{equation}\label{eq:ipa}
    F(\mathbf{x})=H_T\left(H_{T-1}\left(\cdots\left(H_1(\mathbf{x})\right)\right)\right) = H_T\circ H_{T-1}\circ \cdots\circ H_1(\mathbf{x}),
\end{equation}
where $H_1:\Re^d\to\Re^{d_1}$, $H_t:\Re^{d_{t-1}}\to\Re^{d_t}$ for $t=2,\ldots,T-1$, and $H_T:\Re^{d_{T-1}}\to\Re$ are all stochastic functions, and $\circ$ denotes function composition. For simplicity, we can set $d_0=d$ and $d_T=1$. Then, $H_t:\Re^{d_{t-1}}\to\Re^{d_t}$ for all $t=1,\ldots,T$. 

Given the input decision variables $\mathbf{x}$, a simulation algorithm takes $T$ steps. It first calculates $H_1:=H_1(\mathbf{x})$, then $H_2:=H_2(H_1)$, and so on, finally outputting $F(\mathbf{x})=H_T:=H_T(H_{T-1})$. For a function $g:\Re^u\to\Re^v$, it's often reasonable to assume that its computational complexity is in the order of $u\cdot v$. Consequently, the computational complexity of the simulation algorithm can be expressed as
\begin{equation}\label{eq:ipa_sim}
    {\tt Complexity\ of\ Simulation} = \mathcal{O}\left(\sum_{t=1}^T d_{t-1}d_t\right).
\end{equation}

Now consider the IPA algorithm that produces a gradient estimation alongside the simulation process. Let $H_0(\mathbf{x})=\mathbf{x}$ and $D_0$ be the Jacobian matrix of $H_0$ with respect to $\mathbf{x}$. Notice that $D_0$ is a $d\times d$ identity matrix. Let $J_t$ denote the Jacobian matrix of $H_t(H_{t-1})$ with respect to $H_{t-1}$ and let $D_t$ denote the Jacobian matrix of $H_t\circ H_{t-1}\circ\cdots \circ H_1(\mathbf{x})$ with respect to $\mathbf{x}$. It is observed that $J_t$ is a $d_t\times d_{t-1}$ matrix and $D_t$ is a $d_t\times d$ matrix. For each step $t=1,\ldots,T$, once the simulation algorithm computes $H_t$, the IPA algorithm computes $J_t$ and sets $D_t = J_t\cdot D_{t-1}$. Thus, the IPA gradient estimator is $\nabla_{\mathbf{x}} F(\mathbf{x})=D_T'$. Assuming that the effort required to compute an element of any Jacobian matrix, i.e., $\partial H_{t,i}/\partial H_{t-1,j}$, is constant, it becomes evident that the computational complexity of the IPA algorithm is given by
\begin{equation}\label{eq:ipa_ipa}
    {\tt Complexity\ of\ IPA} = \mathcal{O}\left(\sum_{t=1}^T d_{t-1}d_t d\right)=\mathcal{O}\left(d\sum_{t=1}^T d_{t-1}d_t \right).
\end{equation}
Therefore, the computational complexity of the IPA algorithm is an order of $d$ higher than that of the simulation algorithm. In high-dimensional SO, when $d$ is large, the cost of computing IPA estimators may be significantly higher than that of simulation, and this cannot be overlooked.

Another method for computing the sample-path gradient is the back-propagation (BP) algorithm, which is widely used in training neural networks \citep{rumelhart1986learning}. Instead of computing the gradient alongside the simulation process, the BP algorithm calculates the gradient after the entire simulation process is completed. Let $\nabla_{H_T} F=1$ because $F=H_T$. Then, the BP algorithm takes a backward approach, starting from the last step $T$. For each step $t=T-1, T-2,\ldots,0$, the BP algorithm computes $\nabla_{H_t} F = J_{t+1}'\cdot \nabla_{H_{t+1}} F$. It is important to note that $H_0(\mathbf{x})=\mathbf{x}$. Therefore, the BP gradient estimator is $\nabla_{\mathbf{x}} F(\mathbf{x})=\nabla_{H_0} F$. Assuming that the effort of computing an element of any Jacobian matrix, i.e., $\partial H_{t,i}/\partial H_{t-1,j}$, is constant, the computational complexity of the BP algorithm can be expressed as
\begin{equation}\label{eq:ipa_bp}
    {\tt Complexity\ of\ BP} = \mathcal{O}\left(\sum_{t=1}^T d_{t-1}d_t \right).
\end{equation}
Hence, based on (\ref{eq:ipa_sim}), (\ref{eq:ipa_ipa}), and (\ref{eq:ipa_bp}), it can be concluded that the BP algorithm and the simulation algorithm have the same order of computational complexities, while both are an order $d$ faster than the IPA algorithm. Therefore, for high-dimensional problems, the BP algorithm is clearly more advantageous than the IPA algorithm when calculating sample-path gradients.

\cite{wang2023large} consider the more specific example of simulation-based inventory optimization. They demonstrate that the orders of computational complexity for their simulation, IPA, and BP algorithms are $T d^2$, $T d^3$, and $T d^2$, respectively. These correspond to our situation with $d_t=d$ for all $t=1,\ldots,T-1$. In their problems, $d$ may be as high as 500,000. It is not surprising that they observe, compared to the IPA algorithm, the BP algorithm achieves speedups of more than tens of thousands of times when $d$ is large.

\section{Divide and Conquer}\label{DC}

The concept of divide-and-conquer has been effectively employed to tackle large-scale SO problems. Due to the high computational complexity of such problems, divide-and-conquer strategies aim to break down the original problem into smaller sub-problems and address each of these individually. The solutions to these sub-problems are then combined to yield the final solution to the original problem. There are various algorithms that utilize such strategies, differing primarily in terms of the division rules based on problem structures or optimization methods. For instance, in R\&S problems, solutions (or alternatives) may be divided into different groups, while in discrete and continuous SO problems, feasible regions may be partitioned into disjoint sets.

We begin by discussing divide-and-conquer strategies in R\&S procedures in Section \ref{DC_RaS}. This approach involves dividing the solutions into different groups and estimating them in parallel to achieve sample-size efficiency. Additionally, we explore the use of divide-and-conquer strategies in random search and Bayesian optimization algorithms for solving discrete and continuous SO problems. Specifically, we examine the construction of the most promising area for random search algorithms in Section \ref{DC_RS}, and introduce the multiresolution framework for partitioning the entire feasible region into disjoint sub-regions to identify a promising sub-region for Bayesian optimization in Section \ref{DC_SBO}.

\subsection{Divide-and-Conquer in R\&S Algorithms}\label{DC_RaS}

In moderate-scale R\&S problems, the computational cost of comparison and elimination is assumed to be negligible compared to the simulation cost. This is because the number of solutions, denoted as $k$, is usually small, and each simulation replication may take orders of magnitude longer than the $\mathcal{O}(k^2)$ all-pairwise comparisons. However, as $k$ becomes large as in large-scale R\&S problems, the comparison time may become the bottleneck of the R\&S procedures. To address this issue, \cite{ni2017efficient} propose the good selection procedure (GSP), which utilizes the idea of divide-and-conquer. The procedure divides the $k$ solutions into $m$ groups, thereby reducing the computational complexity of comparisons from $\mathcal{O}(k^2)$ to $\mathcal{O}(k^2/m^2)$. Specifically, after dividing into $m$ groups, the intra-group comparison is performed first, and then the local optimal solutions from each group are compared to obtain the global optimal solution. To further improve the elimination efficiency, under the master-worker parallel computing framework, the master retrieves the $m$ local bests from the $m$ groups at the beginning of each local comparison round to find the global best, and then sends the global best to the $m$ groups for additional comparisons.


In addition to improving comparison efficiency, the use of divide-and-conquer in R\&S can also help achieve a lower growth rate in the expected total sample size. \cite{zhong2022knockout} establish a lower bound of $\mathcal{O}(k)$ for the growth rate of the expected sample size for fixed-precision R\&S procedures, which is known as ``sample optimality" and is particularly meaningful for large-scale problems. However, traditional R\&S procedures face challenges in achieving this sample optimality due to their decomposition of the optimization into all-pairwise comparisons, where the solutions compete with each other to eliminate others. To address this sample inefficiency, \cite{zhong2022knockout} introduce the knockout tournament (KT) procedure, which adopts a divide-and-conquer strategy, and proved that the KT procedure achieves the sample optimality. The procedure proceeds in a round-wise manner. In each round, the procedure first pairs the solutions that are still in contention, and then constructs a match between every pair of solutions using an existing R\&S procedure. The winner of each match can advance to the next round of the selection, and the other one is eliminated. Since the procedure follows a decentralized structure, it can be easily implemented in a parallel computing environment. In particular, in a parallel computing environment with $m$ processors, the set of solutions can be equally divided, and the KT procedure can be applied to each subset on a separate processor. After each subset determines a local best solution, a stage-wise procedure is employed to select the global best among them.

In the realm of fixed-budget R\&S procedures, \cite{hong2022solving} introduce the concept of the rate for maintaining correct selection (RMCS), which signifies the minimal growth rate of the total sampling budget required to prevent the probability of correct selection (PCS) from dropping to zero. A lower-order RMCS ensures a relatively high PCS within a limited budget. Thus, RMCS offers a theoretical framework for comparing different fixed-budget procedures in the context of large-scale R\&S. The lower the order of RMCS, the more efficient the procedure. Traditional fixed-budget R\&S procedures, such as the OCBA and its variants, have an RMCS lower bounded by $\mathcal{O}(k\log k)$. In pursuit of better RMCS, \cite{hong2022solving} propose the fixed-budget KT (FBKT) procedure, inspired by the KT procedure, which achieves an improved RMCS to the order of $\mathcal{O}(k)$, known as the optimal rate (see Section \ref{subsec:RS}). To implement the procedure in parallel computing environments with multiple processors, the solutions are equally divided among groups, and each processor is assigned to one group. Similar to the KT procedure, the FBKT procedure conducts a two-phase selection. In the first phase, each processor identifies the local best solution, while in the second phase, processors use the remaining budget to generate additional observations for the local best solutions, selecting the solution with the largest sample mean as the global best one.


\subsection{Divide-and-Conquer in Random Search Algorithms}\label{DC_RS}

Random search algorithms are widely used for solving discrete SO problems. In each iteration, these algorithms typically sample a number of candidate solutions from the neighborhood of the current sample best solution, evaluate them (possibly with some other previously sampled solutions), and select the best solution of this iteration as the current sample best solution. In practice, when dealing with large-scale problems, the focus of algorithm design often lies in exploring promising areas, which are typically constructed based on divide-and-conquer strategies. Many of these algorithms can be proven to have either global convergence or local convergence.


For globally convergent random search algorithms, the nested partitions (NP) framework proposed by \cite{shi2000nested} is well-known. At each iteration, the algorithm divides the feasible region and identifies the most promising area. If a better solution is found in the current most promising area, the area is further explored through partitioning. Otherwise, NP backtracks to its parent region or to the whole feasible region. A similar iterative approach of constructing and partitioning the promising area is employed in the stochastic branch-and-bound algorithm (SB\&B) by \cite{norkin1998optimal, norkin1998branch}. SB\&B is an extension of the branch-and-bound algorithms used in solving deterministic integer optimization problems. The SB\&B algorithm iteratively divides the feasible region into smaller subregions, estimates the boundary of the objective function for these subregions by solving bounding problems, and selects the subregion with the maximum or minimum bound as the promising area. However, finding the estimator of the boundary is challenging, and the partition structure becomes larger, hindering the application of the SB\&B algorithm. To address these issues, \cite{xu2013empirical} propose the empirical SB\&B algorithm. They address the first problem by estimating the boundary based on the performance of the sampled solutions in the subregion. Meanwhile, they adopt the idea of the NP to integrate the regions that are not the most promising area as a region to reduce the large partition structure.

Global convergence in random search algorithms requires the evaluation of all feasible solutions with an infinite number of observations, which is computationally expensive and impractical. In contrast, locally convergent random search algorithms with an implementable stopping rule may be more useful. Unlike globally convergent algorithms, locally convergent algorithms only need to consider the solutions in their local neighborhood to check the local optimality of a solution. An ideal approach is to divide the original feasible region and assign positive probability to a promising subset. For example, the COMPASS algorithm by \cite{hong2006discrete} adopts a unique structure for the most promising area, which includes all solutions closer to $x_{n}^{*}$ than any other simulated solution, where $x_{n}^{*}$ represents the solution with the optimal cumulative sample average at the end of iteration $n$ (i.e., the current sample best solution). COMPASS focuses the search effort in the most promising area, which adapts at each iteration based on the information collected on all simulated solutions. As COMPASS iterates, the most promising area eventually will only contain $x_{n}^{*}$. When this occurs, a statistical local optimality test is performed on $x_{n}^{*}$ and its neighbors to determine whether it is a local optimal solution.


The COMPASS algorithm experiences a significant slowdown as the problem's dimension increases, as discussed in Section \ref{subsec:compass}. To address this challenge, \cite{hong2010speeding} propose replacing uniform sampling of the most promising area with coordinate sampling. This adjustment increases the likelihood of sampling solutions that are close to the current sample best solution. \cite{xu2013adaptive} take a different approach by suggesting the construction of the most promising area using a hyperbox and still uniformly sampling from it. Their algorithm, the adaptive hyperbox algorithm (AHA), defines the most promising area as the largest hyperbox enclosing the current sample best solution and having all other simulated solutions either on the boundary or outside. This construction results in a much smaller hyperbox compared to the most promising area in the COMPASS algorithm, encompassing mostly solutions close to the current sample best. This design enables AHA to scale effectively in high-dimensional problems. Recently, \cite{zhou2023data} apply the AHA algorithm to solve large-scale SO problems in car-sharing service design, demonstrating its effectiveness.


To further enhance the finite-time performance of the COMPASS algorithm and adapt it for commercial solvers, \cite{xu2010industrial} develop an industrial strength COMPASS (ISC) algorithm. This ISC algorithm retains the core of the COMPASS algorithm while incorporating additional steps to enhance its efficiency in solving practical large-scale problems. The optimization process in ISC is divided into three stages: a global search stage, a local search stage, and a final clean-up stage. Specifically, the ISC algorithm utilizes a niching genetic algorithm for the global stage to explore the entire feasible set and identify several promising regions with potentially competitive locally optimal solutions. For the local stage, the COMPASS algorithm is employed to exploit local information and find a locally optimal solution for each of the identified regions from the global stage. Lastly, a R\&S procedure is used for the clean-up stage to select the best solution from all the locally optimal solutions identified. It is noteworthy that within the ISC framework, AHA can also be utilized as the local search algorithm instead of COMPASS.


\subsection{Divide-and-Conquer in Bayesian Optimization}\label{DC_SBO}

Bayesian optimization is a well-established method for addressing discrete and continuous SO problems \citep{sun2014balancing,xie2016bayesian,pearce2022bayesian}. Despite its success in various fields such as material engineering \citep{chen2022sequential}, hyperparameter tuning \citep{snoek2012}, and drug design \citep{negoescu2011}, the application of Bayesian optimization is limited to problems of moderate scale. One of the primary challenges lies in the computational burden associated with estimating the surrogate model \citep{hong2021surrogate}, often a Gaussian process. As the surrogate model involves numerically inverting a large matrix to process the simulation data, it becomes increasingly demanding in computation when the number of design points is large and may eventually become more expensive than the simulation model. To address this challenge, an intriguing idea is to employ a multiresolution framework that adopts a divide-and-conquer strategy. Specifically, the multiresolution framework aims to partition the entire feasible set into disjoint local regions. Subsequently, a promising local region can be identified using a region-level model, and a solution-level model within the selected region is fitted to guide detailed search within this region.


For large-scale discete SO problems, \cite{salemi2019dovs} propose the Gaussian Markov improvement algorithm (GMIA) to identify the global optimal solution. Specifically, GMIA models the objective function at different points as a realization of a discrete Gaussian Markov random field (GMRF) in order to understand the spatial relationships, and then generates a suitable acquisition function known as complete expected improvement (CEI) to balance exploration and exploitation. The GMIA can be applied in a multiresolution framework to enhance the efficiency of the algorithm for problems with large solution spaces. The multiresolution framework operates as follows: it initially divides the solution space into distinct regions and utilizes a region-level GMRF to understand the quality (i.e., response) of these regions. Within each region, it then constructs a solution-level GMRF to learn the quality of all feasible solutions in that region. On one hand, for the region-level GMRF, regions are modeled as nodes in a connected graph, with edges connecting adjacent regions. The response associated with a region (i.e., node) is the average of the objective function values of all solutions in that region, upon which a GMRF can be fitted on these nodes to provide global guidance in finding promising nodes by comparing the corresponding CEIs among all nodes. On the other hand, the quality of individual solutions within a region is represented as a solution-level GMRF, which can help target the most promising solution in this region. This approach significantly reduces the number of points used to calculate the inversion of the covariance matrix for each GMRF, thus saving the computation time. However, such a multiresolution approach may ultimately be limited by the size of the solution-level GMRF it can handle, as it involves calculating CEI for each individual solution. \cite{semelhago2021dovs} propose the rapid Gaussian Markov improvement algorithm (rGMIA) to further enhance solution-level efficiency. This algorithm partitions the solutions in each region into two sets: the search set and the fixed set. In comparison to the original GMIA, rGMIA introduces several rapid search steps. During each rapid search step, it only updates the CEI for the solutions in the search set while fixing the CEI for the solutions in the fixed set. The rationale behind this approach is that the acquisition functions for most solutions remain largely unchanged when new simulation observations are included. By restricting the CEI computation within a small subset of promising solutions for several iterations (i.e., the rapid search steps), computational benefits can be fully exploited.

The multiresolution framework is also well-suited for addressing continuous SO problems. \cite{meng2022combined} introduce the combined global and local search for optimization (CGLO) algorithm, which involves the initial division of the feasible set into separate local regions. A promising region is then selected through optimization over a global model constructed from a finite set of global candidate points. Subsequently, a local search is conducted within the chosen region based on a local model. The algorithm switches back to the global step upon finding a good local solution. In comparison to other Bayesian optimization baseline methods \citep{quan2012simulation,  jalali2017comparison}, the global and local nature of the algorithm enables CGLO to escape suboptimal regions and converge towards the global optimal solution more rapidly, particularly crucial in scenarios with highly multi-modal response functions.

%

\section{Dimension Reduction}\label{dimension reduction}

Solving high-dimensional SO problems, particularly in the pursuit of global optimal solutions, is known to be extremely difficult. This difficulty is largely attributed to the ``curse of dimensionality," which presents itself in several distinct ways. Firstly, as the problem dimension (i.e., the number of decision variables) increases, the search space for optimal solutions expands at an exponential rate. Consequently, SO algorithms must navigate through an increasing number of potential solutions before reaching convergence. Secondly, the number of simulation observations required to obtain a reliable evaluation of the objective function tends to grow exponentially with the dimensionality, often reaching levels that may be unmanageable given available computing power. Finally, SO algorithms may easily become trapped in local optimal solutions, as high-dimensional SO problems are more likely to contain numerous such solutions. To address these challenges, a common strategy involves judiciously transforming a high-dimensional problem into a lower-dimensional one, which is typically easier to solve.


\subsection{Low Effective Dimension}

In many practical scenarios, although the objective function may involve a large number of variables, the actual impact on the function is often exerted by a relatively small subset of these variables. For example, in portfolio optimization, only a few select stocks or bonds may predominantly determine the portfolio's returns or associated risks, despite numerous potential investments being available. Similar patterns can be observed in the field of machine learning, where hyper-parameter optimization for neural networks and deep belief networks frequently reveals that only specific key parameters substantially affect model performance \citep{bergstra2012random}. This phenomenon is referred to as ``low effective dimension" (LED) and forms a strong basis for the application of dimension reduction strategies in solving high-dimensional optimization problems.


One line of research papers exploits the concept of LED by concentrating on identifying the ``effective" variables or, alternatively, screening out the ``ineffective" ones. Effective variables are those that significantly impact the objective function. Once these are determined, optimization can be carried out more efficiently by focusing solely on the effective variables and keeping others at constant or nominal values, resulting in a low-dimensional SO problem. In this context, variable-screening (also known as factor-screening) methods are integrated as a preliminary phase in traditional SO algorithms. These methods typically undertake a series of simulation runs to assess the influence of each input variable on the objective function and statistically eliminate ineffective ones. While various screening methods are available, the ideal one for high-dimensional SO should be able to efficiently identify effective variables from a large pool with minimal simulation runs. Sequential bifurcation (SB), a technique originally applied to simulation by \cite{cheng1997searching} and \cite{Kleijnen2006screening}, enhances efficiency through sequential group screening. It divides the surviving variables into two subgroups at each step and then discards the entire subgroup if it is collectively deemed ineffective. Building upon this, \cite{wan2006controlled} introduce a controlled sequential bifurcation (CSB) procedure, which combines multistage hypothesis testing with the original SB to control both type I error and power for screening. The efficiency of the CSB procedure is further improved in subsequent research by \cite{wan2010improving}.

In parallel, another line of research bypasses the initial step of identifying effective variables and instead focuses on directly streamlining the search process for solutions with LED. Pioneering work by \cite{bergstra2012random} demonstrates that, in the context of hyper-parameter optimization, random search is often more effective than widely used methods such as grid search and manual search. This is attributed to the thorough exploration of the solution space that random search can conduct, is essential for an SO algorithm to locate a global optimal solution. With the assumption of LED, the original high-dimensional space can be embedded into a low-dimensional subspace. This means that dense coverage of the subspace can be achieved through random sampling in each dimension, without the need to determine effective variables in advance. Building on this, \cite{wang2016bayesian} propose a random embedding Bayesian optimization (REMBO) algorithm, which operates within a low-dimensional subspace and projects solutions back to the high-dimensional space. REMBO operates under the assumption that the embedding is linear, though subsequent research \citep[e.g.,][]{Jaquier2020high} has expanded upon this to consider the potential of nonlinear embeddings.


\subsection{Utilizing Special Structures}

In a high-dimensional SO problem where all variables are meaningful, the LED hypothesis might not hold true. Therefore, screening any variable or compressing the problem in a low-dimensional subspace may result in a distorted solution. In such cases, an alternative approach to dimensionality reduction could involve exploiting the inherent structure of the objective function itself.

A common approach is to leverage the additive structure within the objective function, assuming that the high-dimensional objective function can be decomposed into a sum of lower-dimensional functions, each defined on a subset of the underlying variables. \cite{kandasamy2015high} exploit this additive structure by postulating that these subsets are disjoint. This allows each lower-dimensional function to be optimized independently, thus simplifying the high-dimensional optimization challenge. They propose the Additive Gaussian Process Upper Confidence Bound (Add-GP-UCB) algorithm and prove that it can achieve a regret that grows only linearly in the dimension of the original problem. Building on this approach, \cite{rolland2018high} explore a more flexible additive structure that permits overlap between subsets of lower-dimensional components, modeling the interactions between subsets with a graph. However, a challenge in practice is that the decomposition is often unknown and must be inferred. To identify an effective decomposition, \cite{rolland2018high} employ a technique involving random sampling of possible decompositions and selecting the one that maximizes the likelihood. Addressing the vast number of potential additive structures in high-dimensional spaces, \cite{gardner2017discovering} apply a Metropolis–Hastings algorithm to efficiently navigate through the space of possible structures.

Leveraging the smoothness of the objective function for dimension reduction in optimization problems is also a widely adopted strategy. On one hand, we may utilize the gradient and Hessian information to identify low-dimensional subspaces that are critical to optimization. Numerous methods have been developed to explore such information. For instance, \cite{Constantine2014ActiveSubspace} introduce the active subspace methods, which employ the gradient of the objective function to uncover influential directions, enabling the approximation of the objective function within low-dimensional subspaces without sacrificing critical information. \cite{Griebel2006SparseGrids} introduces the sparse grids approach, employing a multi-level hierarchical framework that selectively combines grid points from different resolutions, utilizing mixed derivatives up to the second order to achieve dimension reduction in high-dimensional problems by emphasizing areas of significant variation and de-emphasizing regions of lesser importance. On the other hand, with an understanding of a function's smoothness, we can target a much smaller subspace for efficient approximation. For example, \cite{zhang2021efficient} demonstrate that a function under certain smoothness restrictions lies in the reproducing kernel Hilbert space (RKHS) induced by a Brownian field kernel. The exploration of higher-order derivatives for dimension reduction remains relatively underdeveloped. However, there is an intuitive understanding that as the degree of smoothness in a function increases, the dimensionality of the subspace containing that function decreases (see for detailed discussion of this issue in Section \ref{subsec:smoothness}). Consequently, harnessing insights from higher-order smoothness could provide substantial improvements in dimension reduction, offering a promising avenue for refining optimization strategies.

\section{Gradient-Based Algorithms}

The gradient-based algorithm is a crucial method for continuous simulation optimization. It depends on estimating the gradient of the objective function to determine suitable search directions and utilizes an iterative approach to achieve an optimal solution. The general form of the gradient-based algorithm can be expressed as:
\begin{equation}
\bx_{n+1} = \Pi_{\mathbb{X}} \left( \bx_n - \gamma_n\bK \hat \nabla f(\bx_n) \right),
\end{equation}
where $\nabla f(\bx)$ represents an estimation of the gradient of the objective function $f(\bx)$, $\Pi_{\mathbb{X}}$ indicates a projection back into the feasible region $\mathbb{X}\subseteq \Re^d$, and $\gamma_n$ and $\bK$ are the step size and an appropriate matrix, respectively, used to determine the search length along each dimension of the solution.

Stochastic approximation is a prominent algorithm within the realm of gradient-based algorithms. Initially, it was devised to address the root-finding problem $g(\bx) = 0$ for $\bx\in\mathbb{X}$, and has since found wide-ranging applications in finding (local) optimal solutions or stationary points, such that $\nabla f(\bx) = 0$. When $\hat{\nabla} f$ serves as an unbiased estimator of $\nabla f$, the stochastic approximation algorithm is commonly recognized as being of the Robbins-Monro (RM) type (\citealt{robbins1951stochastic}), where $\{\gamma_n,n=1,2,\ldots\}$ represents a sequence of positive constants satisfying $\sum_{i=1}^{\infty} \gamma_n=\infty \text { and } \sum_{i=1}^{\infty} \gamma_n^2<\infty$. In cases where $\hat{\nabla} f$ is only asymptotically unbiased, for instance, when using a finite difference estimate with the difference approaching zero at an appropriate rate, the algorithm is referred to as being of the Kiefer-Wolfowitz (KW) type (\citealt{KW1952stochastic}). This implies that 
\begin{equation}\label{eq:KWgradient}
\hat{\nabla} f\left(\bx\right)=\left(\begin{array}{c}
\frac{F\left(\bx+c_n\be_1 , \xi_1^{+}\right)-F\left(\bx-c_n\be_1, \xi_1^{-}\right)}{2 c_n} \\
\vdots \\
\frac{F\left(\bx+c_n\be_n , \xi_d^{+}\right)-F\left(\bx-c_n\be_n, \xi_d^{-}\right)}{2 c_n} \\
\end{array}\right),
\end{equation}
where $\be_i$ is a vector with the $i$th component as $1$ and the rest as $0$, and $\xi_i^+$ and $\xi_i^-$ are random variables in the performance function. The sequences $\{\gamma_n,n=1,2,\ldots\}$ and $\{c_n,n=1,2,\ldots\}$ should satisfy 
\begin{equation*}
c_n \rightarrow 0, \quad \sum_{n=1}^{\infty} \gamma_n=\infty, \quad \sum_{n=1}^{\infty} \gamma_n c_n<\infty, \text{ and }\sum_{n=1}^{\infty} \gamma_n^2 c_n^{-2}<\infty.
\end{equation*}

In the realm of large-scale machine learning, stochastic gradient descent (SGD), stemming from stochastic approximation, has emerged as a popular algorithm based on gradients. Consider a data set with $k$ input-output pairs $\{(\bs_i,\by_i),i=1,\ldots,k\}$ and a loss function $\ell(h(\bs,\bx),\by)$, where $h(\bs,\bx)$ and $\by$ denote the predicted and true outputs, respectively. The objective of SGD is to minimize the {\em empirical risk} function $f_n(\bx) = \sum_{i=1}^k \ell(h(\bs_i,\bx),\by_i)/k$. The update rule for SGD is given by
\begin{equation}
    \textbf{x}_{n+1}=\Pi_{\mathbb{X}}(\textbf{x}_n-\gamma_n \hat{\nabla}f_{i_n}(\textbf{x}_n)),
\end{equation}
where $\hat{\nabla}f_{i_n}(\textbf{x}_n)$ represents an estimate of the gradient with respect to the sample $i_n$, and the sample $i_n$ is chosen randomly from $\{1,\cdots,k\}$. Typically, explicit forms of the loss function $\ell$ and the predicted function $h(\bs,\bx)$ are available and, therefore, unbiased estimates of the gradient of the empirical loss can be obtained. The stepsize $\gamma_n$ is the same as in the RM algorithm.

It is evident from the aforementioned algorithms that the estimation of gradients is a critical step in constructing gradient-based algorithms. Therefore, in the subsequent section, we first review the development of gradient estimation methods over the past two decades, followed by a review of gradient-based algorithms for solving large-scale problems.

\subsection{Gradient Estimation}

Gradient estimation, also referred to as sensitivity analysis, is essential not only for developing gradient-based algorithms but also for evaluating the significance of model parameters. In the realm of finance, the gradient of financial product prices in relation to model parameters is sometimes known as Greeks, and is utilized for risk hedging. Therefore, gradient estimation stands as a critical research issue in simulation and has been a subject of study for many years. Gradient estimation methods can be broadly categorized into indirect and direct techniques. Indirect methods primarily involve estimating gradient values using sample functions and gradient representations, with specific methods including the finite difference method and the simultaneous perturbation method. On the other hand, direct gradient estimation methods utilize information from sample functions, their gradients, and the distribution of random variables to estimate gradient values, with perturbation analysis and likelihood ratio (LR)/score function methods being common approaches. For a deeper understanding of these methods, please refer to \cite{fu2015stochastic}.


Over the last two decades, research into gradient estimation has primarily focused on addressing discontinuities in simulation samples, estimating gradients of non-expectation form functions (e.g., quantiles), and developing effective gradient estimators for large-scale problems. In addressing discontinuities, \cite{hong2009pathwise} develop an estimator for probability sensitivities, enabling the estimation of probability sensitivities using the same simulation observations used to estimate probabilities for both terminating and steady-state simulations. Additionally, it applies importance sampling to accelerate the rate of convergence of the estimator. \cite{liu2010kernel} introduce a generalized pathwise method accommodating the discontinuities and propose kernel estimators that require minimal analytical expressions and are easy to implement. \cite{peng2018new} propose a new unbiased gradient estimator, the generalized likelihood ratio (GLR) estimator, to handle discontinuous sample performances with structural parameters. It extends existing methods, such as IPA and LR, to a broader framework. \cite{glynn2021computing} propose a new sensitivity estimator for distortion risk measures that can handle discontinuous sample paths and distortion functions and establish a central limit theorem for the new estimator. \cite{peng2020maximum} develop unbiased estimators for the density and its derivatives for the output of a generic stochastic model using Monte Carlo simulation and propose a gradient-based simulated maximum likelihood estimation method to estimate unknown parameters in stochastic models without assuming an analytical likelihood function.

Previous research on gradient estimation has typically focused on functions in the form of expectation, such as $f(\bx) = \mathds{E}[F(\bx,\xi)]$. There has been relatively less research on functions that are not in the form of expectations, such as quantiles. \cite{hong2009estimating} demonstrates that quantile sensitivities can be expressed as conditional expectations and introduces a batched estimator based on IPA. \cite{fu2009conditional} propose to estimate quantile sensitivities using conditional Monte Carlo simulations and develop a framework to estimate quantile sensitivities efficiently by incorporating conditional expectations and probabilities. \cite{jiang2015technical} provide an alternative derivation of the IPA estimator for quantile sensitivity, simplifying the proofs for strong consistency and convergence rate of the unbatched estimator, and establish strong consistency and a central limit theorem for the batched estimator. Apart from quantiles, \cite{hong2009simulating} derive a closed-form expression for conditional Value-at-Risk (CVaR) sensitivity and introduce an estimator of the CVaR sensitivity with asymptotic properties.


In large-scale optimization problems, such as training neural networks, the BP method is a commonly used approach for gradient estimation. The BP method has a strong connection to the classical IPA method (\citealt{peng2021new}), but its computational complexity surpasses that of IPA, as discussed in Section \ref{subsec:ipa}. Recent works on the BP focus on improving its computational efficiency, such as reducing memory consumption (\citealt{gruslys2016memory}) and meta learning (\citealt{NEURIPS2021_7608de7a}). Apart from the BP, automatic differentiation, initially proposed by  \cite{linnainmaa1970representation}, has also garnered attention and is now integrated into many programming languages and packages (\citealt{van2017tangent, paszke2017automatic}). Automatic differentiation can be seen as an interpretation of a computer program that includes the calculation of derivatives alongside the standard computation process. Essentially, any numerical computation can be broken down into a series of basic operations for which derivatives are well defined. By applying the chain rule to combine the derivatives of these basic operations, we can determine the derivative of the entire computation (\citealt{rall1996introduction}).

\subsection{Gradient-Based Algorithms in Large-Scale Problems}

Despite being an early proposed method, gradient-based algorithms have undergone significant advancements in both methodology and application over the past two decades. Particularly with the rapid progress in computer science, gradient-based algorithms have become an important method for solving large-scale optimization problems such as training neural networks. Therefore, in this subsection, we focus on reviewing the development of gradient-based algorithms, especially those of the SGD, in addressing large-scale problems. For further exploration into the development of other variants of gradient-based algorithms, please refer to \cite{lan2020} for more details.

Generally, gradient descent (GD) is a commonly used method for training machine learning models. However, for problems with large data sizes, SGD proves superior to GD. \cite{bottou2010large} reveals the tradeoffs for the case of small-scale and large-scale learning problems and demonstrates the superiority of SGD in large-scale problems. Specifically, let $\varepsilon >0$ denote the preset training error. The total number of iterations required to achieve an optimization accuracy of $\varepsilon$ is proportional to $\log(1/\varepsilon)$ when the loss function is strongly convex. This implies that, with $k$ samples, the total computational cost for achieving $\varepsilon$-optimality for GD is proportional to $k\log(1/\varepsilon)$. Conversely, for SGD, the total computational cost of $\varepsilon$-optimality is proportional to $1/\varepsilon$. Notice that the $\varepsilon$-optimality of SGD does not depend on the size of the data set $k$, so SGD is a more suitable approach for handling large-scale datasets compared to GD.

Recent research has primarily focused on two aspects to enhance the efficiency of SGD. The first aspect aims to alleviate the negative impact of noisy gradient estimates, which can hinder SGD's convergence under a fixed stepsize setting and lead to slow, sub-linear convergence rates when using decreasing stepsizes. Previous studies have explored three noise reduction methods to address this issue: dynamic sampling, gradient aggregation, and iterate averaging. Dynamic sampling techniques progressively enlarge the mini-batch size during gradient computation to reduce noise, leading to more precise gradient estimates as the optimization progresses (\citealt{Byrd2012,HashemiGP2014}). Gradient aggregation techniques enhance the accuracy of the search directions by retaining gradient estimates corresponding to samples utilized in prior iterations, revising some of these estimates in each iteration, and defining the search direction as a weighted average of these estimates. Representative variants include stochastic variance reduced gradient \citep{Johnson2013}, stochastic average gradient (\citealt{RouxSB2012}), and SAGA (\citealt{DefazioBL2014}). Iterate averaging techniques achieve noise reduction by maintaining an average of iterates computed during the optimization process, rather than averaging gradient estimates (\citealt{Polyak1990,polyak1992,NemirovskiJLS2009}).


The second aspect involves mitigating the adverse effects of the high nonlinearity and ill-conditioning of the objective function by integrating second-order information. These techniques enhance the convergence rates of batch methods or improve the constants in front of the sub-linear convergence rate of the SGD. Algorithms incorporating second-order information typically take the following form:
\begin{equation}
        \textbf{x}_{n+1}=\Pi_{\mathbb{X}}(\textbf{x}_n-\gamma_n \textbf{K}_k^{-1} \hat{\nabla}f_{i_n}(\textbf{x}_n)).
\end{equation}
Here, $\textbf{K}_k$ represents the second-order information of the objective function. Intuitively, methods such as the Gauss-Newton method and the quasi-Newton method compute the Hessian matrix or its inverse and use it as the second-order information. The Gauss-Newton method generates an estimate of the Hessian using only first-order data, ensuring positive semi-definiteness even in cases where the Hessian is indefinite (\citealt{schraudolph2002fast}). Quasi-Newton methods, like BFGS (Dennis and Schnabel, 1983) and L-BFGS (\citealt{dennis1974characterization}) and L-BFGS (\citealt{liu1989limited,schraudolph2007stochastic}), dynamically update a Hessian approximation using an algorithm instead of computing the actual Hessian at each iteration. These methods are effective for optimizing problems with millions of variables, making them suitable for large-scale optimization problems. In addition to Hessian-based second-order information, the natural gradient method defines a search direction in the space of realizable distributions rather than within the space of the real parameter vector, utilizing the Fisher information matrix as the second-order information (\citealt{park2000adaptive,marceau2016practical,martens2020new}). Furthermore, to tackle larger-scale problems, some algorithms aim to avoid matrix multiplications. For instance, the conjugate gradient method updates iterations using Hessian vector products instead of the full Hessian matrix (\citealt{dembo1982inexact,byrd2011use,roosta2019sub}), while RMSprop, Adadelta, and AdaGrad methods restrict attention to diagonal or block-diagonal scaling matrices (\citealt{duchi2011adaptive,tieleman2012rmsprop,zeiler2012adadelta,kingma2014adam}).

\section{Parallelization}

In practical terms, large-scale SO problems require simulating thousands to millions of solutions, each with multiple observations, which may result in unmanageable computational costs. With the advancement of multi- to many-core environments, it is natural to harness parallel computing to expedite large-scale SO. There are two distinct approaches to parallel SO. The first approach involves parallelizing the optimization process while maintaining the simulation of a single replication on a single processor. The second approach encompasses parallelizing not only the optimization process, but also the simulation itself. For the first approach, the challenges lie in synchronizing different simulation replications handled on different processors. These challenges are discussed for parallel R\&S in Section \ref{PRS}, and parallel Bayesian optimization in Section \ref{PBO}, respectively. As for the second approach, the difficulty lies in finding a general parallel computing framework that is applicable to different simulation models. The concept of tensorization is discussed to take advantage of widely available parallel computing tools and to significantly reduce the barrier of parallelization in Section \ref{tensorization}. It is important to note that parallelization naturally incorporates the idea of divide-and-conquer, which has been extensively introduced in Section \ref{DC}. Therefore, this section focuses more on the problems created by parallelization and how to address them, rather than on the implementation of parallelization in the divide-and-conquer scheme.

\subsection{Parallel R\&S}\label{PRS}

R\&S procedures typically involve simulating each alternative with multiple replications and then comparing the average performance of the alternatives to determine the best option. It is clear that these procedures are well-suited for parallelization due to the ease of parallelizing the simulation replications. However, difficulties arise when the replication times of different alternatives are different or stochastic, making it challenging to synchronize multiple processors. Careless handling of this issue can lead to significant loss of efficiency.

Many R\&S procedures attempt to address this issue by avoiding synchronization. In a master-worker parallel computing environment, widely used for parallel R\&S procedures, a master processor dispatches the jobs (i.e., simulation replications for different alternatives) to multiple worker processors. The asynchronous handling of the simulation observations typically means that the order of the job completions (referred to as the output sequence) from the workers may differ from the order in which the master assigns the jobs (referred to as the input sequence). This misalignment may lead to significant statistical issues for classical sequential R\&S procedures when the simulation run times depend on the simulation results \citep{luo2015fully}. Several studies in the literature have focused on adapting the existing fully-sequential procedures (e.g., the KN procedure of \cite{kim2001fully}) to accommodate this situation. One natural approach is to restore the order of the input sequence in the output sequence. This approach is termed ``vector filling" by \cite{luo2015fully} and ``zipping" by \cite{ni2013ranking}. However, this approach may require a large amount of memory and could potentially lead to issues in the event of communication interruptions (i.e., significant delays in returning a simulation result to the master). To tackle these challenges, \cite{luo2015fully} introduce the asymptotic parallel selection (APS) procedure, which demonstrates that the error resulting from ignoring the difference between the input and output sequences is negligible in a reasonable asymptotic regime.

In addition to the synchronization issue, all-pairwise comparisons (often needed for fully-sequential R\&S procedures) on the master and the communications between the master and workers are also crucial concerns in parallel R\&S \citep{hong2021review}, especially for very large-scale R\&S problems. \cite{zhong2022speeding} propose the parallel Paulson's procedure (PPP), which modifies the well-known Paulson's procedure \citep{paulson1964sequential}, and breaks all-pairwise comparisons into comparisons with the ``best", thus reducing the complexity of comparison calculation from $\mathcal{O}(k^2)$ to $\mathcal{O}(k)$, which is different from and more efficient than the divide-and-conquer method used by \cite{ni2017efficient} (also introduced in Section \ref{DC_RS}). To reduce the frequent communication between processors, the PPP procedure divides the surviving alternatives into $m$ groups and instructs each worker to perform the simulations for one group of alternatives, significantly reducing the communication times.

Similarly, \cite{pei2018new} introduce a framework known as parallel adaptive survivor selection (PASS), which involves individual comparisons of alternatives to a standard that dynamically adapts as surviving alternative observations accumulate. In the bisection-PASS (bi-PASS) variant \citep{pei2022parallel}, the standard is estimated by the weighted sample mean of current not-yet-eliminated alternatives. Within such a framework, the complexity of comparison becomes $\mathcal{O}(k)$, and the elimination process becomes more efficient than other procedures, due to the large number of eliminated alternatives at each iteration. It is evident that both the PASS and bi-PASS procedures facilitate relatively fast communication, as only a few scalar quantities are communicated between the master and any worker. However, these procedures do not guarantee a single selection, but aim to retain all good alternatives. To assess the performance of the retained set, \cite{pei2022parallel} propose a new guarantee on the expected false elimination rate (EFER), which emphasizes retaining a proportion of good systems relative to the total number of good systems, instead of using the classical PCS guarantee.

\subsection{Parallel Bayesian Optimization}\label{PBO}

Traditional Bayesian optimization algorithms typically comprise two primary components: a surrogate model and an acquisition function. The surrogate model, often implemented as a Gaussian process, offers a Bayesian posterior probability distribution of the objective function $f(\mathbf{x})$ at a candidate solution $\mathbf{x}$. Upon observing $f$ at a new solution, this posterior distribution is updated. The acquisition function gauges the potential effectiveness of a solution $\mathbf{x}$ if it were to be evaluated, based on the current posterior distribution of $f$. Subsequently, the acquisition function is optimized to select the most suitable solution for evaluation next.

Bayesian optimization algorithms are inherently sequential, typically evaluating one solution at a time. When implementing these algorithms in parallel computing environments, it becomes necessary to evaluate multiple solutions simultaneously. While this does not pose a significant challenge for building surrogate models, it does complicate the selection of suitable solutions using the acquisition function. Two fundamental approaches exist to address this issue: the synchronous approach and the asynchronous approach.

The synchronous approach involves selecting a batch of solutions based on the acquisition function, assigning them to multiple processors, and waiting until all of them are evaluated before moving to the next iteration. While the concept of this approach is straightforward, the challenge lies in optimizing the acquisition function to select multiple solutions. Many methods propose addressing this challenge by considering a multi-solution approximation of the original acquisition function that is easier to optimize, or by proposing approximation algorithms that solve the optimization problem approximately. For instance, \cite{ginsbourger2008} extend the expected improvement (EI) function to the $q$-EI fucntion, allowing the selection of $q$ solutions simultaneously. \cite{wang2020parallel} develop an efficient algorithm, combining the IPA gradient estimator and multi-start gradient descent, to optimize the $q$-EI function. \cite{shah2015parallel} extend the acquisition function of \cite{hernandez2014predictive} and propose the predictive entropy search (PES). \cite{wu2016parallel} develop the parallel version of the knowledge gradient algorithm (referred  as $q$-KG) and proposed to maximize the $q$-KG using a similar algorithm as \cite{wang2020parallel}.

The asynchronous approach does not wait for all solutions to be evaluated. Instead, it selects the next solution to evaluate whenever a processor becomes idle, significantly reducing waiting time, especially when the simulation times of different solutions are different or stochastic. However, the challenge lies in considering the solutions being evaluated on other processors. Ignoring these solutions may result in the selection of a solution similar to the running solutions, thereby reducing the efficiency of parallel computing. To overcome this challenge, some studies attempt to promote diversity in exploration by integrating the location information of the running solutions into the acquisition function. For instance, building on the Gaussian process upper confidence bound (GP-UCB) algorithm proposed by \cite{srinivas2010gaussian}, \cite{desautels2014parallelizing} introduce the Gaussian process batch upper confidence bound (GP-BUCB) algorithm, which takes into account the impact of the running solutions on the posterior variance of the Gaussian process. By updating the variance using all the solutions, including those whose observations are not yet available, GP-BUCB can encourage exploration in regions that have not been sampled yet. \cite{wang2020parallel} constrain the remaining solutions in the $q$-EI optimization to be different from the previous simulated solutions as well as the running solution, thus ensuring some variation among them.

\subsection{Tensorization of Simulation Model}\label{tensorization}

The previously mentioned approaches concentrate on parallelizing the optimization process, but do not address parallelizing the individual simulation replications. However, there are simulation optimization challenges where the simulation itself poses a difficulty. This scenario often arises when the simulation is utilized to model an exceedingly complex system, and its execution may consume a considerable amount of time using traditional simulation algorithms. For example, \cite{wang2023large} examine the inventory optimization issue in a production system with potentially 500,000 products, where traditional simulation algorithms may require hours to execute a single replication. With the evolution of deep learning, various computational tools like TensorFlow \citep{tensorflow2015-whitepaper} and PyTorch \citep{paszke2019pytorch} have been designed for rapid tensor operations. These tools can effectively harness the parallel computing capabilities offered by multi-core CPUs and many-core GPUs to expedite calculations in complex models. Therefore, to achieve substantial acceleration for large-scale simulation models, a crucial lesson gleaned from the success of deep learning is the tensorization of the simulation model. In this section, we explore the recent advancements in the parallelization of single replications across diverse simulation models using tensorization.

The inventory simulation model discussed above shares a critical resemblance with recurrent neural networks (RNN), as noted by \cite{wang2023large}. They propose an RNN-inspired simulation approach that significantly accelerates simulations, enabling the resolution of large-scale inventory optimization problems within reasonable timeframes. A key aspect of their approach involves tensorization of the simulation model, enabling the use of efficient computational tools available in software packages and substantially reducing the barriers to parallelization. Their numerical experiments demonstrate that while tensorization does not enhance the theoretical computational complexity, it dramatically reduces the runtime of the simulation algorithm by orders of magnitude.

The simulation of queueing networks often relies on discrete-event simulation as a popular method. However, parallelizing queue network simulation using discrete-event simulation presents challenges due to the need for a global event list to maintain the correct event order and capture the effects of event interactions. Addressing the effective management of large-scale queueing networks, \cite{hong2024fast} introduce a simulation method tailored for large-scale Markovian queueing networks by incorporating an Euler approximation. This method can accommodate time-varying dynamics and is optimized for efficient implementation through vectorization (or tensorization) techniques. By formulating the algorithms in terms of vector and matrix operations, they enable the utilization of parallel computing capabilities inherent in multi-core CPU and many-core GPU architectures, along with leveraging efficient computational tools designed for vector operations. The use of vectorization allows for the simultaneous handling of multiple node tasks, leading to a significant acceleration in computation speed and rendering this method highly suitable for simulating large-scale queueing networks.

\section{Concluding Remarks}

There are several remarks that we want to make to conclude this review article. Firstly, large-scale SO problems differ fundamentally from their moderate counterparts, requiring different theories, algorithms, and, above all, a completely different mindset to solve them. Secondly, large-scale R\&S problems have been extensively researched, with potential for further enhancing the efficiency of algorithms. However, the potential for fundamental breakthroughs may be limited. Thirdly, insufficient attention has been given to high-dimensional SO, including both discrete and continuous problems, leaving ample room for the development of new theories and algorithms. We believe that both dimension reduction and the utilization of high-order smoothness are crucial for the development of efficient high-dimensional SO algorithms. Finally, parallelization is essential for large-scale SO, and the effective utilization of many-core GPUs presents a significant opportunity for advancing theory and algorithm design in this field.

\section*{Statement}
During the preparation of this paper, the authors used ChatGPT in order to improve the readability and the English language of this paper. After using this service, the authors reviewed and edited the content as needed and take full responsibility for the content of the paper.

\bibliographystyle{apalike}
\bibliography{ref}

\end{document}